\documentstyle[11pt,a4,amscd,amssymb,verbatim]{amsart}
\newtheorem{Thm}{Theorem}
\newtheorem{Rem}{Remark}
\newtheorem{Cor}{Corollary}
\newtheorem{lemma}{Lemma}

\newcommand{\Ima}{\operatorname{Im}}

\newcommand{\Hom}{\operatorname{Hom}}
\newcommand{\End}{\operatorname{End}}

\newcommand{\modc}{\operatorname{-mod}}

\newcommand{\Span}{\operatorname{span}}
\newcommand{\seq}{\operatorname{seq}}

\newcommand{\res}{\operatorname{res}}
\newcommand{\TL}{\operatorname{TL}}
\newcommand{\C}{\mathbb C}

\edef\savecatcodeat{\the\catcode`@}
\catcode`\@=11

\def\tb@ifSpecChars#1#2{#1}
\def\tb@ifNoSpecChars#1#2{#2}

\def\tableau{%
  \bgroup% matched in \tb@tableauD
  \@ifstar{\let\Tif\tb@ifNoSpecChars\tb@tableauB}% *, don't use special chars
          {\let\Tif\tb@ifSpecChars\tb@tableauB}}% no *, use special chars

\def\tb@tableauB{% add [] if no [options]
  \@ifnextchar[{\tb@tableauC}{\tb@tableauC[]}}

\def\tb@tableauC[#1]{\hbox\bgroup%
    \let\\=\cr% end line
    \def\bl{\global\let\tbcellF\tb@cellNF}%
    \def\tf{\global\let\tbcellF\tb@cellH}% highlighted cell
    \dimen2=\ht\strutbox \advance\dimen2 by\dp\strutbox%
    \ifx\baselinestretch\undefined\relax%
    \else%
       \dimen0=100sp \dimen0=\baselinestretch\dimen0%
       \dimen2=100\dimen2 \divide\dimen2 by\dimen0%
    \fi%
    \let\tpos\tb@vcenter% default position
    \tb@initYoung% default tableau type
    \tb@options#1\eoo% parse options
    \let\arrow\tb@arrow%
    \dimen0=\Tscale\dimen2%
    \dimen1=\dimen0 \advance\dimen1 by \tb@fframe%
    \lineskip=0pt\baselineskip=0pt% line spacing will be from \vbox to \dimen0
    \def\tb@nothing{}%
    \def\endcellno{$\rss\egroup\bss\egroup}% end cell w/o overlap
    \def\endcell{\endcellno\kern-\dimen0}% end cell & prepare to overlap it
    \def\begincell{\vbox to\dimen0\bgroup\vss\hbox to\dimen0\bgroup\hss$}%
    \let\overlay\tb@overlay%
    \let\fl\tb@fl%
    \let\lss\hss\let\rss\hss\let\tss\vss\let\bss\vss% cell alignment
    \def\mkcell##1{% format individual cell
        \let\tbcellF\tb@cellD% default cell frame
        \def\tb@cellarg{##1}% store cell contents
        \ifx\tb@cellarg\tb@nothing\let\tb@cellarg\tb@cellE\fi%
        \begincell\tb@cellarg\endcellno% the actual cell content
        \tbcellF}% draw cell frame
    \let\savecellF\tbcellF% save global value of cellF in case of nested tableau
     \Tif{\catcode`,=4\catcode`|=\active}{}\tb@tableauD}%

\let\tb@savetableauD\tableauD% save any current definition
{% set up characters which will be interpreted as command characters
    \catcode`|=\active \catcode`*=\active \catcode`~=\active%
    \catcode`@=\active% command characters
\gdef\tableauD#1{%
  \Tif{% make all the command characters active in math mode when #1 parsed
    \mathcode`|="8000 \mathcode`*="8000%
    \mathcode`~="8000 \mathcode`@="8000%
    \def@{\bullet}%
    \let|\cr% end line
    \let*\tf% highlighted cell
    \let~\sk% skew cell
  }{}%
  \tpos{\tabskip=0pt\halign{&\mkcell{##}\cr#1\crcr}}%
  \global\let\tbcellF\savecellF% restore global value
  \egroup% match \hbox\bgroup at start of \tableauC
  \egroup}% match \bgroup at start of \tableau
}
\let\tb@tableauD\tableauD% rename the command
\let\tableauD\tb@savetableauD% restore old command with this name
\let\tb@savetableauD\undefined

\def\tb@options#1{\ifx#1\eoo\relax\else\tb@option#1\expandafter\tb@options\fi}

\def\tb@option#1{%
  \if#1t\let\tpos\tb@vtop\fi%        t = align at top
  \if#1c\let\tpos\tb@vcenter\fi%     c = align at center
  \if#1b\let\tpos\vbox\fi%           b = align at bottom
  \if#1F\tb@initFerrers\fi%          F = Ferrers diagram
  \if#1Y\tb@initYoung\fi%            Y = Young diagram
  \if#1s\tb@initSmall\fi%            s = small boxes
  \if#1m\tb@initMedium\fi%           m = medium boxes
  \if#1l\tb@initLarge\fi%            l = large boxes
  \if#1p\tb@initPartition\fi%            p = small partition sized boxes
  \if#1a\tb@initArrow\fi%            a = use arrow font as base dimension
}

\def\tb@vcenter#1{\ifmmode\vcenter{#1}\else$\vcenter{#1}$\fi}

\def\tb@vtop#1{\hbox{\raise\ht\strutbox\hbox{\lower\dimen0\vtop{#1}}}}

\def\tb@initPartition{\def\Tscale{.3}}
\def\tb@initSmall{\def\Tscale{1}}
\def\tb@initMedium{\def\Tscale{2}}
\def\tb@initLarge{\def\Tscale{3}}

\def\tb@initArrow{\dimen2=1.25em}

\def\tb@initYoung{%
  \def\tb@cellE{}% empty cell stays empty
  \let\tb@cellD\tb@cellN% default frame is normal frame
  \def\sk{\global\let\tbcellF\tb@cellNF}}% skew cells are empty
\def\tb@initFerrers{%
  \def\tb@cellE{\bullet}% empty cell gets bullet
  \let\tb@cellD\tb@cellNF% default frame is no frame
  \def\sk{\bullet}}% skew cell gets bullet

\tb@initMedium% default scale

\def\tb@sframe#1{%
  \vbox to0pt{%            Embed frame in a box of no vert or hor extent
    \vss%                            pull box above reference point
    \hbox to0pt{%
      \hss%                          pull box left of reference point
      \vbox to\dimen1{%              Actual width of frame
        \hrule depth #1 height0pt% draw top edge of frame
        \vss%                     begin vcenter sides
        \hbox to\dimen1{%           horiz box with side edges just inside
          \vrule width #1 height\dimen1% left edge
          \hss%                     stretch center
          \vrule width #1%         right edge
          }%
        \vss%                     end vcenter sides
        \hrule height #1 depth 0in% bottom edge
        }%
      \kern-\tb@hframe%           horiz alignment off by half line width
      }%
    \kern-\tb@hframe}}%           vert alignment off by half line width
\def\tb@hframe{.2pt}\def\tb@fframe{.4pt}\def\tb@bframe{1.2pt}
\def\tb@cellH{\tb@sframe{\tb@bframe}}       % bold frame
\def\tb@cellNF{}                            % no frame
\def\tb@cellN{\tb@sframe{\tb@fframe}}       % normal frame
\let\tbcellF\tb@cellN                       % default is normal

\def\tb@rpad{1pt}
\def\tb@lpad{1pt}
\def\tb@tpad{1.8pt}
\def\tb@bpad{1.8pt}

\def\tb@overlay{\endcell\@ifnextchar[{\tb@overlaya}{\begincell}}
\def\tb@overlaya[#1]{\vbox to\dimen0\bgroup%
  \tb@overlayoptions#1\eoo%
  \tss\hbox to\dimen0\bgroup\lss}
\def\tb@overlayoptions#1{\ifx#1\eoo\relax\else\tb@overlayoption#1\expandafter\tb@overlayoptions\fi}

\def\tb@overlayoption#1{
  \if#1t\def\tss{\vskip\tb@tpad}\let\bss\vss\fi% t = align at top
  \if#1c\let\tss\vss\let\bss\vss\fi%             c = align at center
  \if#1b\def\bss{\vskip\tb@bpad}\let\tss\vss\fi% b = align at bottom
  \if#1l\def\lss{\hskip\tb@lpad}\let\rss\hss\fi% l = align at left
  \if#1m\let\lss\hss\let\rss\hss\fi%             m = align at middle
  \if#1r\def\rss{\hskip\tb@rpad}\let\lss\hss\fi% r = align at right
}

\def\tb@fl{\endcell\begincell\vrule depth 0pt width \dimen0 height \dimen0 \endcell\begincell}

\@ifundefined{diagram}{}{
\def\tb@arrowpad{.5}

\newoptcommand{\tb@arrow}{\@ne}[2]{%
  \endcell% end previous cell contents
   \begingroup%
   \let\dg@getnodesize\tb@getnodesize% substitute routine to get nodesize
   \dg@USERSIZE=#1\relax%
   \ifnum\dg@USERSIZE<\@ne \dg@USERSIZE=\@ne \fi%
   \dg@parse{#2}%
   \dg@label{\tb@draw{#1}{#2}}}% draw arrow

\def\tb@getnodesize#1#2#3#4#5{\dimen3=\tb@arrowpad\dimen2 #4=\dimen3 #5=\dimen3\relax}
\def\tb@getnodesize#1#2#3#4#5{\ifnum#2=0\ifnum#3=0\tb@getnodesizetail{#4}{#5}\else\tb@getnodesizehead{#4}{#5}\fi\else\tb@getnodesizehead{#4}{#5}\fi}
\def\tb@getnodesizetail#1#2{\dimen3=.5\dimen2 #1=\dimen3 #2=\dimen3}
\def\tb@getnodesizehead#1#2{\dimen3=.5\dimen2 #1=\dimen3 #2=\dimen3}

\def\tb@draw#1#2#3#4{%
        \dg@X=0\dg@Y=0\dg@XGRID=1\dg@YGRID=1\unitlength=.001\dimen0%
        \dg@LBLOFF=\dgLABELOFFSET \divide\dg@LBLOFF\unitlength%
        \dg@drawcalc% compute arrow geometry
        \begincell% start tableau cell
        \let\lams@arrow\tb@lams@arrow% substitute routine
        \begin{picture}(0,0)\begingroup\dg@draw{#1}{#2}{#3}{#4}\end{picture}%
        \endcell% end tableau cell
        \endgroup% match \begingroup in \tb@arrow
        \begincell}% start new entry in this cell
}

\def\tb@lams@arrow#1#2{%
 \lams@firstx\z@\lams@firsty\z@
 \lams@lastx#1\relax\lams@lasty#2\relax
 \lams@center\z@
 \N@false\E@false\H@false\V@false
 \ifdim\lams@lastx>\z@\E@true\fi
 \ifdim\lams@lastx=\z@\V@true\fi
 \ifdim\lams@lasty>\z@\N@true\fi
 \ifdim\lams@lasty=\z@\H@true\fi
 \NESW@false
 \ifN@\ifE@\NESW@true\fi\else\ifE@\else\NESW@true\fi\fi
 \ifH@\else\ifV@\else
  \lams@slope
  \ifnum\lams@tani>\lams@tanii
   \lams@ht\ten@\p@\lams@wd\ten@\p@
   \multiply\lams@wd\lams@tanii\divide\lams@wd\lams@tani
  \else
   \lams@wd\ten@\p@\lams@ht\ten@\p@
   \divide\lams@ht\lams@tanii\multiply\lams@ht\lams@tani
  \fi
 \fi\fi
 \ifH@  \lams@harrow
 \else\ifV@ \lams@varrow
 \else \lams@darrow
 \fi\fi
}

\catcode`\@=\savecatcodeat
\let\savecatcodeat\undefined
\begin{document}
\large

\title{The Ariki-Terasoma-Yamada tensor space and the blob algebra. } 
\author{Steen Ryom-Hansen}
\footnote{Supported in part by Programa Reticulados y Ecuaciones and by FONDECYT grants 1051024
and 1090701.} 
\address{Instituto de Matem\'atica y F\'isica, Universidad de Talca \\
Chile\\ steen@@inst-mat.utalca.cl }
\maketitle

\begin{abstract}
We show that the Ariki-Terasoma-Yamada tensor module 
and its permutation submodules $ M(\lambda) $ are modules for the blob algebra 
when the Ariki-Koike algebra is a Hecke algebra of type $B$.
We show that $ M(\lambda)$ and the standard modules $ \Delta(\lambda) $
have the same dimensions, the same localization and 
similar restriction properties and are equal in the Grothendieck 
group. Still we find that the universal property for 
$ \Delta(\lambda) $ fails for $ M(\lambda) $, making $ M(\lambda) $ and $ \Delta(\lambda) $
different modules in general.
Finally, we prove that $ M(\lambda) $ is isomorphic to the dual Specht module for the 
Ariki-Koike algebra.

\end{abstract}

\section {\bf Introduction. }
In this paper we combine the representation theories of the Ariki-Koike algebra and 
of the blob-algebra. 
The link between the two theories is 
the tensor space module $ V^{\otimes n} $ for the Ariki-Koike algebra 
defined in [ATY] by Ariki, Terasoma and Yamada.

\medskip

The blob algebra $ b_n = b_n(q,m) $ was defined by Martin and Saleur [MS] as a 
generalization of the Temperley-Lieb algebra
by introducing periodicity in the statistical mechanics model. 
The blob algebra is 
also sometimes called the Temperley-Lieb algebra of type $ B $, or the
one-boundary Temperley-Lieb algebra, 
and indeed it has a diagram calculus 
generalizing the Temperley-Lieb diagram calculus. Our work 
treats the non-semisimple representation theory of $ b_n$.

\medskip
There is a natural embedding $ b_n \subset b_{n+1}   $ which 
gives rise to restriction and induction functors between the module categories. 
These functors are part of a powerful  
category theoretical formalism on the representation theory of the entire 
tower of algebras. It also 
involves certain localization and globalization functors $ F $ and $G $ between the 
categories
of $b_n$-modules for different $n$. We denote it the localization/globalization formalism. 

\medskip

The formalism is closely related to the fact that $b_n $ 
is quasi-hereditary in the sense of Cline, 
Parshall and Scott, [CPS], (when $ q + q^{-1} \not= 0 $).
Its parametrizing poset is 
$ \Lambda_n := \{ n , n - 2, \ldots  , -n \} $.
The standard modules $ \Delta_n(\lambda), \, \lambda \in \Lambda_n $ 
can be defined by a diagram basis and 
have dimensions equal to certain binomial coefficients. 

\medskip
A main point of our work is the existence of 
a surjection $ \pi $ from the Hecke algebra $ H(n,2) = H_n(q, \lambda_1, \lambda_2 )$ 
of type $B_n$ to the blob algebra $ b_n$, for appropriate choices of the parameters. 
It makes it possible to pullback $ b_n $-modules to $ H(n,2) $-modules and in this way 
the category of $ b_n$-modules may be viewed as a subcategory of 
the $ H(n,2) $-modules. 

\medskip

Since $ H(n,2) $ is a special case of an Ariki-Koike algebra 
it has a tensor module $ V^{\otimes n } $ as described in [ATY].
As a first result we prove that $ V^{\otimes n } $ and its 
'permutation' submodules $ M_n(\lambda) $
are $ b_n $-module when $ \dim V = 2 $.
We are then in position to apply the 
localization/globalization formalism to the module $ M_n(\lambda) $, 
and to compare it to the
standard module $ \Delta_n(\lambda) $. 

\medskip
In our main results we show that the two modules have the same dimensions, 
share the same localization properties and even are equal 
in the Grothendieck group of $ b_n $-modules. They also have 
related behaviors under restriction from $ b_n $ to $ b_{n-1} $. 
Even so we find that $ M_n(\lambda) $ and $ \Delta_n(\lambda) $ are different modules
in general. 
We show this by demonstrating that the universal 
property for $ \Delta_n(\lambda) $ fails for $ M_n(\lambda) $. To be 
more precise, we show that in general
$ G F M_n(\lambda) \not\cong M_n(\lambda) $
whereas it is known that $ G F \Delta_n(\lambda) \cong \Delta_n(\lambda) $ 
(when $ \lambda \not= \pm n $).

\medskip
This rises the question whether 
$ M_n(\lambda)$ may be identied with another 
'known' module. 
We settle this question by considering the Specht module 
$ S(n_1, n_2)$ 
for  
$ H(n,2)$, 
where $ (n_1, n_2 ) $ is a two-line bipartition associated with $ \lambda $.
We show that this module is the pullback of a $ b_n $-module, also 
denoted $ S(n_1, n_2)$, and that $ M_n(\lambda)$ is isomorphic to the contragredient dual 
of $ S(n_1, n_2)$.

\medskip
We find that, somewhat surprisingly, neither of the $ b_n$-modules
$ M_n( \lambda) $, 
$  S(n_1, n_2) $ nor	 
their duals identify with the standard module
$ \Delta_n(\lambda) $ for $ b_n$.
 
\medskip
It is pleasure to thank P. Martin for useful conversations. 
Thanks are also due to the referee for useful comments. 

\section{\bf Preliminaries}

In this section we shall briefly recall the results of [MW] and [ATY], the two main sources of 
inspiration for the present paper.
Let us start out by the work of Martin-Woodcock [MW]. Among other things they realize the blob 
algebra $b_n $
as a quotient of the Ariki-Koike algebra $ H(n,2) $ by the ideal generated by the idempotents associated with certain 
irreducible representations of $ H(2,2) $. It then turns out that this ideal has 
a simple description in terms of the $ H(n,2) $-generators. Let us explain all this briefly.

\medskip

Let $ {\cal A }= {\mathbb Z }[q,q^{-1}, \lambda_1, \lambda_2] $. 
Let $ H(n,2)= H(n,q, \lambda_1 , \lambda_2 ) $ be the unital
$ \cal A $-algebra generated by 
$ \{ X, g_1, \ldots , g_{n-1} \} $ with relations 
$$ 
g_i g_{i \pm 1 }g_i = g_{i \pm 1 } g_i g_{i \pm 1 } \,\,\,\,\,\,\,\,\,\,\,\,\, [ g_i, g_j ] = 0 \,\,\,\, i \not= j \pm 1 
$$
$$
g_1 X g_1 X = X g_1 X g_1   \,\,\,\,\,\,\,\,\,\,\,\,\,\,\,\,\, [ X , g_j ] = 0 \,\, \, \, j > 1  
$$
$$
(g_i - q ) (g_i + q^{-1} )=0 
$$
$$
(X-\lambda_1)(X-\lambda_2) = 0 $$
It is the $ d=2 $ case of the Ariki-Koike algebra $  H (n,d) $ 
or the cyclotomic Hecke algebra of type $ G(d,1,n)$, see 
[AK] and [BM]. 
For $ \lambda_1 = - \lambda_2^{-1} $ it is the Hecke algebra of type $ B_n$. 
Note that there is a canonical embedding $ H(n, 2) \subset H(n+1, 2) $.

\medskip
As usual, if $k$ is an $ \cal A $-algebra we write $ H_k(n,2) := H(n,2) \otimes_{\cal A} k$
for the specialized algebra.

\medskip
Recall the concept of cellular algebras, that was introduced by Graham and Lehrer in [GL]
in order to provide a common framework for many algebras that appear in non-semisimple 
representation theory. It is shown in [GL] that the Ariki-Koike algebra 
is cellular for 
general parameters $ n,d$.
In our case $ d= 2 $ it also follows from [DJM].

\medskip
Let $ k $ be a field and suppose that $ k $ is made into an $ \cal A $-algebra by mapping 
$ q, \lambda_1, \lambda_2 $ to nonzero elements $ q, \lambda_1, \lambda_2 $ of $ k $.
Assume that $ q^4 \not= 1 $, 
$ \lambda_1 \not= \lambda_2 $ and $  \lambda_1 \not= q^2 \lambda_2 $. Then there are formulas for
$ e^{- 1}, e^{- 2} \in {  H}_k(2,2) $, the primitive central idempotents corresponding to the two one-dimensional cell 
representations given by $ ( 1^2, \emptyset ) $, $ (  \emptyset, 1^2 ) $, see [MW] for a more precise statement concerning the 
actual cell modules that we are refering to and for the details.
Let $ I \subset H_k(n,2) $ be the ideal in $ { H}_k(n,2) $
generated by $ e^{- 1}, e^{- 2} $.
%and define $ H^{\cal D}(n,2) $ as the quotient $  H (n,2)/ I $. 
Using the mentioned formulas, it is shown in (27) of [MW] that $ I $ is generated by either of
the elements 
$$ 
\begin{array}{c}
(X_1+X_2 -(\lambda_1 + \lambda_2)) (g_1-q)  \\
(X_1X_2 -\lambda_1 \lambda_2) (g_1-q) 
\end{array}
$$
where as usual $ X_1 := X, \, X_i := g_{i-1} X_{i-1} g_{i-1} $ for $ i = 2,3 \ldots $.

\medskip

Let $ m \in \mathbb Z $ and assume that $ n $ is a positive integer.
The blob algebra $ b_n = b_n(q,m) $ is the unital $ k $-algebra 
on generators $ \{ U_0, U_1 , \ldots U_{n-1} \} $ and 
relations $$  U_i U_{i \pm 1 } U_i = U_i, \, U_i^2 = -[2]\, U_i, \, U_0^2 = -[m]\, U_0, \, U_1 U_0 U_1 = [m-1] U_1  $$
for $ i > 0 $ and commutativity between the generators otherwise. As usual $ [ m ] $ is here the Gaussian integer 
$ [m] :=  \frac{q^m - q^{-m}}{q - q^{-1} }$. The blob algebra was introduced in [MS] via a basis of decorated Temperley-Lieb algebras, which 
explains it name. We shall however mostly need the above presentation of it. 
This is only one of several different presentations of $ b_n $, the one used in [MW].

\medskip
Let $ H^{\cal D}(n,2) $ be the quotient $  H_k(n,2)/ I $ 
and choose $$ \lambda_1 = \frac{q^m}{q-q^{-1}} \,\,\,\,\,\, \mbox{and} \,\,\,\,\,\,
 \lambda_2 = \frac{q^{-m}}{q-q^{-1}} $$ 
Using the above description of $ I $, 
it is then shown in proposition (4.4) of [MW] that 
the map $ \varphi $ given by 
$ \varphi:  g_i -q \mapsto  U_i, \, \, \,   X - \lambda_1 \mapsto   U_0 $ 
induces a $k$-algebra isomorphism 
\begin{equation}{\label{MW-iso}}
 \varphi :H^{\cal D}(n,2) \cong b_n(q,m)
\end{equation}

We finish this section by recalling the construction of the tensor representation of the Ariki-Koike algebra $  H (n,d) $
found by Ariki-Terasoma-Yamada [ATY]. It is an extension to the Ariki-Koike case of Jimbo's 
classical tensor representation 
of the Hecke algebra, [J],
and therefore basically amounts to the extra action of $ X $ factorizing 
through the relations. On the other hand, this action is quite 
non-trivial and is for example not 
local in the sense of [MW].

\medskip

The [ATY] construction works for all Ariki-Koike algebras $ {H} (n,d) $, but we shall only
need the $ d= 2 $ case, which we now explain.
Let 
$ V $ be a free $ \cal A$-module of rank two and let $ v_1, v_2 $ be a basis.
Let 
$ R \in \End_{\cal A} (V \otimes V) $ be given by 
$$  \left\{ \begin{array}{l} R( v_i \otimes v_j ) =   
q v_i \otimes v_j  \,\,\,\,\,\,\,\,\,\,   \mbox{if } i=j  \\
R( v_2 \otimes v_1 ) = 
v_1 \otimes v_2 \\ R( v_1 \otimes v_2 ) = 
v_2 \otimes v_1 + (q-q^{-1})  v_1 \otimes v_2 
\end{array} \right\}
$$
Then the $ H(n,2) $ generator $ g_i $ acts on $ V^{\otimes n } $ 
through $$  T_{i+1} :=Id^{\otimes i-1 } \otimes R \otimes Id^{\otimes n-i-1} $$ 
The $ g_i $ generate 
a subalgebra of $ H(n,d) $ isomorphic to the Iwahori-Hecke algebra of type $ A $ and the above action 
is the $ \dim V = 2 $ case of 
the one found by Jimbo in [J]. The maximal quotient of it acting faithfully 
on $ V^{\otimes n} $ is the Temperley-Lieb algebra $ \TL_n $.

\medskip
For $ j =2, 3 , \ldots, n $ we shall need the 
$ \cal A$-linear map
$ S_j \in \End_{{\cal A}} (V^{ \otimes n } ) $, that 
by definition acts on 
$ v = v_{i_1} \otimes v_{i_2} 
\otimes \cdots \otimes v_{i_{j-1}} \otimes v_{i_{j}} \otimes \cdots \otimes v_{i_n } $ 
through 
$$ S_j( v ) = \left\{
\begin{array}{ll} q v   
& \mbox{if   }  i_{j-1} = i_j  \\
v_{i_1} \otimes v_{i_2} \otimes \cdots \otimes v_{i_{j}} \otimes v_{i_{j-1}} \otimes \cdots 
\otimes v_{i_n}   & \mbox{otherwise }
\end{array} \right.
$$
Let $ \theta := S_n S_{n-1} \cdots S_2 $ 
and let 
$ \varpi \in \End_{ {\cal A}}	(V^{ \otimes n } ) $ be the map given by
$$ 
v_{i_1} \otimes v_{i_2} 
\otimes \cdots \otimes v_{i_n } 
\mapsto 
\lambda_{\delta(1)}
v_{i_1} \otimes v_{i_2} 
\otimes \cdots \otimes v_{i_n } 
$$
where $ \delta(1) = 1 $ if $ i_1 = 1 $ and 
$ \delta(1) = 2 $ if $ i_1 = 2 $.	
Then 
$ \theta \varpi $ is given by  
$$ \theta \varpi: \, \,  v_{i_1} \otimes v_{i_2} \otimes v_{i_3} \otimes \cdots \otimes v_{i_n} \mapsto \lambda_{\delta(1)}
q^{a-1}  v_{i_2} \otimes v_{i_3} \otimes \cdots \otimes v_{i_n} \otimes v_{i_1} $$ 
where $ a $ is the number of $ i_k $ such that $ i_k = i_1 $. 
Now [ATY] define the action of $ X  \in H(n,2)$ by the formula 
$$ T_1 := T_2^{-1} T_3^{-1} \cdots T_n^{-1} \theta \varpi $$
As mentioned in 
[ATY], the proof that the $ T_1, T_2, \ldots, T_{n-1} $ satisfy the Ariki-Koike relations 
works in specializations as well. 
One of the steps of their proof 
is the following lemma, which we shall need later on.
\medskip

\begin{lemma}~\label{basiclemma}
Let $ Y_{j,p} $  be the $ \cal A $-submodule of $ V^{\otimes n } $ generated by basis elements 
$ v= v_{i_1} \otimes v_{i_2} \otimes \cdots \otimes v_{i_n} $ such that $ i_p \geq j $. Then
if $ v \in Y_{j,p} $ we have that 
$$ T_{p+1}^{-1} T_{p+2}^{-1} \cdots T_n^{-1} 
S_{n} S_{n-1} \cdots S_{p+1} \, v  = v  \mbox{ mod } \,  Y_{j+1,p} $$
\end{lemma}
\section{\bf The Ariki-Terasoma-Yamada tensor space as blob algebra module}

From now on we assume that $ k $ is an algebraically closed field, 
such that $ q, \lambda_1, \lambda_2 \in k $ and $q^4 \not= 1 $, 
$ \lambda_1 \not= \lambda_2,  \lambda_1 \not= q^2 \lambda_2 $. 
We moreover assume that $ \lambda_1 = \frac{q^{m}}{q-q^{-1}} $ and
$ \lambda_2 = \frac{q^{-m}}{q-q^{-1}} $ where $ m $ is an integer.
With these assumptions the results of the previous section are valid.

\medskip

In this section we prove 
that the Ariki-Koike action given by the above construction factors through the
blob algebra. Let $ V $, $ T_i $ be as in the previous section. Then we have 
\begin{Thm}{\label{ATYblob}}
$ (T_1 T_2 T_1 T_2 - \lambda_1 \lambda_2 )(T_2 -q ) = 0 $ in $ \End_k(V^{\otimes n }) $.
\end{Thm}
\begin{pf*}{\it Proof}
We start by noting that by the Ariki-Koike relations 
$$ (T_1 T_2 T_1 T_2 - \lambda_1 \lambda_2 )(T_2 -q ) = (T_2 -q ) (T_1 T_2 T_1 T_2 - \lambda_1 \lambda_2 )   $$
We show that $ (T_1 T_2 T_1 T_2 - \lambda_1 \lambda_2 )(T_2 -q )=0  $ 
on all basis elements of $ V^{\otimes n } $. 
It clearly holds for $ v = v_{i_1} \otimes v_{i_2} \otimes \cdots \otimes v_{i_n} $ where $ i_1 = i_2 $, so 
we assume $ i_1 \not= i_2 $. If $ i_1 = 2 $ and $ i_2 = 1 $ we get by lemma \ref{basiclemma}
that the action of $ T_1 $ on 
$ v $ is multiplication by $ \lambda_2 $. But then $ T_2 T_1 T_2 $ acts on $ v $ through 
$$ \begin{array}{c} T_2 T_1 T_2 
(v_{ 2 } \otimes v_{ 1 } \otimes \cdots \otimes v_{i_n}) = \\
T_{2}  T_{2}^{-1} T_{3}^{-1} \cdots T_n^{-1} 
S_{n} S_{n-1} \cdots S_{2} \varpi T_{2} 
(v_{2} \otimes v_{1} \otimes \cdots  \otimes v_{i_n} )=  \\
\lambda_1 T_{3}^{-1} \cdots T_n^{-1} 
S_{n} S_{n-1} \cdots S_{3} 
(v_{2} \otimes v_{1} \otimes \cdots \otimes v_{i_n} ) = \\
\lambda_1 (v_{2} \otimes v_{1} \otimes \cdots \otimes v_{i_n} ) \mbox{  mod  } Y_{2,2} 
\end{array}
$$
by lemma \ref{basiclemma} once again. Actually, since $ T_{3}^{-1} \cdots T_n^{-1} 
S_{n} S_{n-1} \cdots S_{3} $ does not change the first coordinate of $ v $ we can even calculate 
modulo the subspace $ Y_2 $ of $ V^{\otimes n } $ generated by $ v_{2} \otimes v_{2} \otimes v_{i_3} \otimes \cdots 
\otimes v_{i_n}  $.
We conclude that $  (T_1 T_2 T_1 T_2 - \lambda_1 \lambda_2 )v  \in  Y_2 $. 
But clearly $ T_2 -q$ kills $ Y_2 $ and we are done
in this case.

\medskip
On the other hand, we have that 
$$ V^{\otimes n } =  \ker(T_2 -q) + \Span_k \{ v_2 \otimes v_1 \otimes v_{i_3} \otimes \cdots \otimes v_{i_n} \, | \,
i_j =1,2  \mbox{ for } j \geq 3 \} $$ 
and hence $ V^{\otimes n } $ is also equal to 
$$  \ker(T_2 -q)(T_1 T_2 T_1 T_2 - \lambda_1 \lambda_2)  + 
\Span_k \{ v_2 \otimes v_1 \otimes v_{i_3} \otimes \cdots \otimes v_{i_n} \, | \,
i_j =1,2  \mbox{ for } j \geq 3 \} 
$$
Combining with the above, the theorem follows.
\end{pf*}
\begin{Rem}
The formula of the theorem is easy to implement in a computer system and amusing to verify.
\end{Rem}

\begin{Cor}
$ V^{\otimes n} $ is a 
$b_n(q,m) $-module with 
$ U_i, \, i \geq 1 $ acting 
through  $  T_{i+1} -q $ and $  U_0 $ through  $  T_1 - \lambda_1 $. 
\end{Cor}
\begin{pf*}{\it Proof}
Using that $ \lambda_1 = \frac{q^{m}}{q-q^{-1}} $ and
$ \lambda_2 = \frac{q^{-m}}{q-q^{-1}} $ (and the other assumptions on the parameters)
this follow from the theorem and proposition (4.4) of [MW].
\end{pf*}

\section{ \bf Localization and globalization}
The main results of our paper depend on a category theoretical approach to 
the representation theory of $ b_n $ that we shall now briefly explain.
It was introduced by J. A. Green in 
the Schur algebra setting, [G], but has turned out to be useful in the context of  
diagram algebras as well, see e.g. [CVM], [M] and [MR]. In the case of the blob algebra $ b_n $,
a good references to the formalism is [MW1], see also [CGM]. 

\medskip
Recall first that $ [2] \not= 0 $ in $ k $ so that we can 
define $ e := -\frac{1}{[2] }U_{n-1} $. This is an idempotent of $ b_n$ and we have 
that $ e b_n  e \cong b_{n-2} $, see [MW1]. Hence it gives rise to the exact localization functor $$ F: b_n \modc \rightarrow b_{n-2} \modc, \, \, \, 
M \mapsto e M $$
It has a left adjoint, the globalization functor 
$$ G: b_{n-2} \modc \rightarrow b_{n}  \modc, \, \, \, 
M \mapsto b_n e \otimes_{e b_n e}  M $$
which is right exact.
Let $ \Lambda_n := \{ n, n-2, \ldots , -n+2, -n \} $.
Under our assumption $ [2] \not= 0 $, the category $ b_n$-mod is quasi-hereditary with labeling poset $ (\Lambda_n, \prec) $, where
$ \lambda \prec \mu \Leftrightarrow | \lambda | > |\mu | $. Hence for all $ \lambda \in \Lambda_n $ we have a standard module 
$ \Delta_n(\lambda) $, a costandard module $ \nabla_n(\lambda) $, a simple module $ L_n(\lambda) $, a projective module 
$ P_n(\lambda) $ and an injective module $ I_n(\lambda) $. The simple 
module $ L_n(\lambda) $ is the unique simple quotient of $ \Delta_n(\lambda) $. In general 
$ \Delta_n(\lambda) $ and $ L_n(\lambda) $ are different.

\medskip
One can find in [MW1] a diagrammatical 
description of $ \Delta_n(\lambda) $. 
We shall however first of all need the following category theoretical 
properties
of $ \Delta_n(\lambda) $. 
Assume first that $ n \geq 3$ to avoid $ b_n$ for $ n \leq 0 $ that 
we have not defined. Then we have 
\begin{equation}{\label{categorial}}
\begin{array}{l}
F \Delta_n(\lambda) \cong \left\{ \begin{array}{ll} 
\Delta_{n-2}(\lambda) & \mbox{if} \,\,\,   \lambda \in \Lambda_{n}  \setminus \{\pm n\}\\
0  & \mbox{otherwise}  
\end{array}
\right. \\
G \circ F \Delta_n(\lambda) \cong \left\{ \begin{array}{ll} 
\Delta_{n}(\lambda) & \mbox{if} \,\,\,   \lambda \in \Lambda_{n}  \setminus \{\pm n\}\\
0  & \mbox{otherwise}  
\end{array}
\right.
\end{array}
\end{equation}
where the second isomorphism is the adjointness map of the pair $ F $ and $ G$. 
Note that the second statement is false if $ \Delta_n(\lambda) $ is
replaced by $ \nabla_n(\lambda) $.
Together with 
$$ \Delta_{n}(\pm n ) \cong L_n(\pm n)  \cong \nabla_n(\pm n) $$
and 
$$
\begin{array}{l}
F L_n(\mu) \cong \left\{ \begin{array}{ll} 
L_{n-2}(\mu) & \mbox{if} \,\,\,   \mu \in \Lambda_{n}  \setminus \{\pm n\}\\
0  & \mbox{otherwise}  
\end{array}
\right. 
\end{array}
$$
these properties give the universal property for $ \Delta_n(\lambda) $. For assume that $ N $ is a $ b_n$-module 
with $ [N: L_n(\lambda)] = 1 $ satisfying $ [N: L_n(\mu)] \not= 0 $ only if $ \mu \prec \lambda $. 
Then applying a sequence of functors $ F $ until arriving at $ L_{| \lambda |}(\lambda ) $ followed by a 
similar sequence of functors $ G $, we obtain a nonzero homomorphism $ \Delta_n(\lambda) \rightarrow N $. In other 
words, $ \Delta_n(\lambda) $ is projective in the category of $b_n$-modules whose simple factors are 
all of the form $ L_n(\mu) $ with $\mu \preceq \lambda $.

\medskip

Let us now return to the tensor space module $ V^{\otimes n } $ 
for $b_n $ from the previous section.
For $ \lambda \in 
\Lambda_n $, we denote by $ M(\lambda) = M_n(\lambda) $ the 'permutation' module. 
By definition, its 
basis vectors are $ v_{i_1} \otimes v_{i_2} \otimes \cdots 
\otimes v_{i_n} $ satisfying 
$$ \lambda = \# \{ j \, | \, i_j =1 \} - \# \{ j \, | \, i_j =2 \} $$ 
It is clear from the previous section that it is a $ b_n $-submodule of $ V^{\otimes n } $.

\medskip

We shall frequently make use of 
the {\it sequence} notation that was introduced in [MR] for the basis vectors of $ V^{\otimes n }$. 
Under it $ 112 $ corresponds to 
$ v_1 \otimes v_1 \otimes v_2 $ and so on. 
As in [MR] the set of 
sequences of $ 1 $s and $ 2 $s of length $ n $ is denoted $ \seq_n $. 
The subset of these sequences with 
$ 1 $ appearing $ n_1 $ times is denoted 
$ \seq_n^{n_1}  $. 
With this notation $ M_n (\lambda) $ has basis $ \seq_n^{ a}  $ 
where $ a = \frac{ \lambda +n}{2} $.
Its dimension 
is given by the binomial coefficient 
$ \left( \begin{array}{c} n \\  a  \end{array} \right) $.
This is also the dimension of $\Delta_n(\lambda)$. 

\medskip

We shall also need 
the underline notation from [MR].
It is useful for doing 
calculations in $ F M  $ where $ M $ is a submodule of $V^{ \otimes n} $.
In the present setup it is given by 	
$ \underline{12} := q^{-1} 12 -21 $ for $ n = 2 $ and extended linearly to 
higher $ n$.  
For example,  
for $ n = 3 , \lambda = 1 $ we get the following identities 
in $ F M_n( \lambda) = e M_n( \lambda) $
$$ 1\underline{12} = [2] e(112) = -U_2(112) = -(T_3-q)(112) = -(121 - q^{-1} 112) $$ 

\medskip
Since $ M_n (\lambda) $ and $\Delta_n(\lambda)$ have the same dimension 
one might guess that they are isomorphic $b_n$-modules. To see whether this is true one 
would have to verify for $ M_n (\lambda) $ the category theoretical properties given
in ({\ref{categorial}}). 
The following theorem shows that the first of these indeed holds.

\medskip
\begin{Thm}{\label{firststep}} For $ n \geq 3 $ 
there is an isomorphism of $ b_{n-2} $-modules
$$ F M_n(\lambda) \cong \left\{ \begin{array}{ll}
M_{n-2}(\lambda) & \mbox{if} \,\,\,   \lambda \in \Lambda_{n}  \setminus \{\pm n \}\\
0  & \mbox{otherwise}  
\end{array}
\right. $$
\end{Thm}
\begin{pf*}{\it Proof} 
The theorem is easy to verify for $ \lambda = \pm n $ so let us assume that 
$ \lambda \in \Lambda_{n}  \setminus \{\pm n \} $.
Let 
$ f: M_{n-2}(\lambda)  \rightarrow F M_n(\lambda)   $ be the $ k$-linear map given by  
$$ i_1 i_2 \cdots  i_{n-2} \mapsto 
i_1 i_2 \cdots  i_{n-2} \underline{12} := 
q^{-1} \, i_1 i_2 \cdots  i_{n-2} 12 -
i_1 i_2 \cdots  i_{n-2} 21
$$
We show that $ f $ is a $ b_{n-2} $-linear isomorphism.

\medskip

But by lemma 1 of [MR] 
we already know that $ f$ is a vector space isomorphism and that it is linear
with  
respect to the Temperley-Lieb action.
Hence we must show that $ f $ is linear with respect to the action of $ X $. 
Here $ X $ acts on the left hand side through the restriction to 
$ M_{n-2}(\lambda) $ of 
$ T_1 \in \End_k ( V^{ \otimes n -2 } ) $ whereas
it acts on the right hand side through the restriction to $ F M_{n}(\lambda) $ 
of $ \frac{-1 \,\, \, }{[2]} (T_n-q) \, T_1 \, \frac{-1 \, \,\,}{[2]}(T_n -q) \in \End_k ( V^{ \otimes n } ) $. Since we assume $ n \geq 3 $ the factors of the product commute. Noting furthermore that 
$ \frac{-1\,\,\,}{\, [2]}(T_n -q) $ acts through the identity on 
$ F M_{n}(\lambda) $, we get that the action of $ X $ on the right hand side is nothing but 
the restriction of 
$ T_1  \in \End_k ( V^{ \otimes n } ) $ to $ F M_{n}(\lambda) $.

\medskip
It is now enough to show that $ f $ is linear with respect to 
$ T_1  \in \End_k ( V^{ \otimes n-2} ) $ and $ T_1  \in \End_k ( V^{ \otimes n } ) $, 
in other words that 
$$ f( T_2^{-1}  \cdots T_{n-2}^{-1} S_{n-2} \cdots S_2 \varpi  v )  
= 
 T_2^{-1}  \cdots  T_{n-1}^{-1} T_{n}^{-1} S_n S_{n-1}  \cdots S_2 \varpi f(v )  $$
for all $ v \in M_{n-2}( \lambda) $. 
For this we first note that $ f $ clearly commutes 
with 
$ T_2, \cdots, T_{n-2},  S_2, \cdots, S_{n-2},$ and $ \varpi$. 
Since these are all invertible, we are reduced to proving that 
\begin{equation}{\label{n=3}}
f(v) = T_{n-1}^{-1} T_{n}^{-1} S_n S_{n-1} f(v) \, \, \mbox{ for all } \, v \in M_{n-2}(\lambda) 
\end{equation}
This equation only involves the last three factors of $ f(v) $ so we may assume that $ n=3$.
But for $ n=3 $, 
the cases $ \lambda =  \pm  3 $ of (\ref{n=3})
are trivially fulfilled, leaving us 
the $ \lambda =  \pm 1 $ cases.

\medskip

If $ \lambda = 1 $
we have that 
$$ \Ima f =  e M_3(1) = \Span_k \{\,1 \underline{12} \,\}  = \Span_k \{\,112-q121 \,\} 
$$
and we must prove that $ T_2^{-1} T_3^{-1} S_3 S_2 (112 - q121 ) = 112 - q121 $ or 
\begin{equation}\label{equivalently}
S_3 S_2 (112 - q121 ) = T_3  T_2 (112 - q121) 
\end{equation}
The left hand side of this equation is $ q(121 - q 211 ) $ whereas the right hand side is 
$$
\begin{array}{cr} 
T_3  T_2 (112 - q121) = T_3 ( q 112 - q ( 211 + (q-q^{-1}) 121 )) & = \\
T_3 ( q 112 - q 211 -  (q^2-1) 121 ) & = \\  q 121 + (q^2-1)112 - q^2 211 -  (q^2-1) 112 =
q 121  - q^2 211    & 
\end{array}
$$
as claimed.

\medskip
If $ \lambda = -1 $ we have that 
$$ \Ima f =  e M_3(-1) = \Span_k \{\,212-q221 \,\} $$ 
and so the equivalent of equation (\ref{equivalently})
is 
$$ S_3 S_2 (212 - q221 ) = T_3  T_2 (212 - q221) $$
The left hand side of this is $ q ( 112 - q 212) $, and the right hand side is 
$$
T_3  T_2 (212 - q221) = T_3 (  122 - q^2 \,221)    = 
  q 122 - q^2 \, 212 
$$
as claimed. The theorem is proved.
\end{pf*}

We now go on to consider the analogue for $ M_n(\lambda)$ of 
the second category theoretical property for $ \Delta_n( \lambda ) $ in (\ref{categorial}). 
It turns out {\it not} to hold for $ M_n(\lambda)$. Let us be more 
precise. Let $ \seq_n^{n_1} $ be the basis for $ M_n(\lambda)$ as above and
define $ n_2 := n - n_1 $ such that 
$ \lambda = n_1 - n_2 $.
We then have the following result.
\begin{lemma}{\label{relate}} 
Let $ n \geq 3 $ and suppose that $ q $ is an $l$th primitive root of unity, where
$ l $ is odd.
Suppose $ \lambda \in \Lambda_{n} \setminus \{\pm n \} $. Then we have 
\newline
\noindent
a) The adjointness map $ \varphi_{\lambda}: G \circ F  M_n(\lambda) \rightarrow M_n(\lambda) $ 
is surjective if and only if $ n_2 \not= m \mod l $. \newline
b) The adjointness map $ \varphi_{\lambda} : G \circ F M_n (\lambda) \rightarrow M_n(\lambda) $ 
is injective iff $ n_2 \not= m \mod l $. \newline
c) The adjointness map $ \varphi_{\lambda} : G \circ F M_n (\lambda) \rightarrow M_n(\lambda) $ 
is an isomorphism iff $ n_2 \not= m \mod l $. 
\end{lemma}
\begin{pf*}{\it Proof}
Part c) obviously follows by combining a) and b). Let us now 
prove a). Assume first that $ n_2 \not= m \mbox{\,mod} \, l $ and suppose that $ \varphi_{\lambda} $ 
is not surjective. 

\medskip
Note first that 
for $ w \in \seq_{n-2}^{n_1 -1} $ and 
$ (i_{n-1}, i_{n}) = (1,2 ) $ or $ (2,1) $   
we have that $ e (w i_{n-1} i_{n})  =   c \, w \underline{12} $ for some scalar $ c \in k^{\times} $.
Recall next from [MR] that $ b_n e $ is generated as an $ e b_n e $ 
right module by the set $$ {\cal G} := \{ U_{n-1},\, U_{n-2} U_{n-1}, 
\,  \ldots, 
U_{0} \cdots U_{n-2}  U_{n-1} \} $$ and that $ \varphi_{\lambda}: G \circ F M_n(\lambda) 
\rightarrow M_n(\lambda) $ is the multiplication map 
$$ b_n e \otimes_{e b_n e} e M_n(\lambda)  \rightarrow M_n(\lambda), 
\, \, \, \, \, \, \, \, \, \, \, \, 
U \otimes m \mapsto Um $$
Suppose that $ w = i_1 i_2 \cdots i_{n-2} $. 
A key point, used in [MR] as well, is now that for $ j \geq 1 $ the multiplication 
of $ U_j U_{j+1} \cdots U_{n-1} \in {\cal G} $ on 
$ w \underline{12} $ shifts the underline 
to position $ (j, j+1) $ in the following sense
$$ U_j U_{j+1} \cdots U_{n-1} w \underline{12} = -[2] \, i_1 i_2 \cdots i_{j-1} \underline{12} 
i_{j+2} \cdots i_n $$
as follows easily from the definitions.
Using it 
we get that $ im  \varphi_{\lambda} $ is the span of $$ I_1 = \{ (X-\lambda_1) \underline{12}x 
\, | \, x \in \seq_{n-2}^{n_1 -1} \, \} $$ together with 
$$ I_2=  \{ v_1  \underline{12} v_2 \, | \, v_1 \in \seq_k^{l_1  }, \, v_2 \in \seq_{n-2-k}^{n_1 - l_1 -1},   \, k \leq  n-2, \, \,  l_1  \leq n_1 -1 \} $$

\medskip
Let $ N_2 := \Span_k \{ w \, | \, w \in I_2 \} $. 
Then $ Q:= M_n(\lambda )/ N_2 $ is a vector space of dimension one since the
 elements of $ I_2 $ can be viewed as
straightening rules 
that allow us to rewrite 
any element of $ M_n(\lambda )/ N_2 $ as a 
scalar multiple of $ 1^{n_1 } 2^{n_2} $ (or $ 2^{n_2 } 1^{n_1} $). Indeed, 
by the definition of $ \underline{12} $ we have the 
following identity, valid in $ Q $
\begin{equation}{\label{straightening}}
 v_1 12 v_2 = q v_1 21 v_2 \, \, \mbox{ for } 
v_1 \in \seq_k^{l_1  }, \, v_2 \in \seq_{n-2-k}^{n_1 - l_1 -1}
\end{equation}
But $ N_2  \subseteq im \varphi_{\lambda}  $ and so we conclude $ im \varphi_{\lambda} = N_2 $ 
since $ \varphi_{\lambda}$ is not surjective.

\medskip

%\medskip

But then $ Q$ is 
a $ b_n $-module. It has dimension one and hence the action of $X$ on $ Q$ is given by a scalar, which we shall work out.
Notice first that if
$ i \geq 2 $ then 
$ T_i^{-1}$ acts through the constant $ q^{-1} $ on $ Q$, since $ U_i $ acts as zero for $ i > 0 $.

\medskip

Set $ v = 1^{n_1} 2^{n_2}  \in Q   $.
Since $ X $ acts through $ T_2^{-1} T_3^{-1} \cdots T_n^{-1} \theta \varpi $ we get 
that 
$$
\begin{array}{c}
X v = \lambda_1 q^{n_1 -1} q^{-n_1 -n_2 +1} 1^{n_1 -1} 2^{n_2}1 = \\
\lambda_1 q^{n_1 -1} q^{-n_1 -n_2 +1} q^{-n_2} 1^{n_1 } 2^{n_2} = 
\lambda_1 q^{-2 n_2 } 1^{n_1 } 2^{n_2} = \lambda_1 q^{-2 n_2 } v  
\end{array}
$$
using the straightening rules (\ref{straightening}).
Hence the scalar in question is $ \lambda_1 q^{-2n_2 } $.

\medskip

Set now $ w = 2^{n_2} 1^{n_1}  \in Q   $. Then we get the same way
$$
\begin{array}{c}
X w = \lambda_2 q^{n_2 -1} q^{-n_1 -n_2 +1}  2^{n_2-1}   1^{n_1} 2   = \\
\lambda_2 q^{n_2 -1} q^{-n_1 -n_2 +1} q^{n_1 } 2^{n_2} 1^{n_1}      =                                              
 \lambda_2  w
\end{array}
$$
The two scalars must be same, that is $ \lambda_1 q^{-2n_2 } = \lambda_2 $ and hence $ \lambda_1 / \lambda_2 = q^{2m} = q^{2n_2}$. 
Since $ l $ is odd, this implies that $ n_2 = m \mod l $, which is the desired contradiction.

\medskip
To prove the other implication we assume that $ n_2 = m \mbox{ mod } l $ 
and must show that $ \varphi_{\lambda} $ is not surjective. We show that 
$ I_1 \subseteq N $ or equivalently $ (N_1 + N_2)/N_2  =0 $ where $ N_1 := \Span_k \{ w \, |\,  w \in I_1 \} $.

\medskip
Since the actions of $ X $ and $ U_i $ commute for $ i = 3, \ldots ,n-1 $, we get for any 
$ w \in \seq_{n-2}^{n_1 -1} $ that 
$$ (X-\lambda_1) \underline{12}w =c X \underline{12} 1^{n_1-1} 2^{n_2 -1} \mbox{ mod } N_2 $$
where $ c \in k^{\times} $.
We go on calculating modulo $ N_2 $ and find 
$$ \begin{array}{rl} 
X \underline{12} 1^{n_1-1} 2^{n_2 -1} = & X  q^{-1} 12  1^{n_1-1} 2^{n_2 -1} - X 21  1^{n_1-1} 2^{n_2 -1}   \\
= & q^{-n_2 -1 } \lambda_1 2  1^{n_1-1} 2^{n_2 -1}1  - \lambda_2   q^{-n_1  } 1^{n_1} 2^{n_2 }   \\
= & q^{-2 n_2 -n_1 } \lambda_1 1^{n_1} 2^{n_2 }  - \lambda_2   q^{-n_1  } 1^{n_1} 2^{n_2 }  = 0  \\
\end{array} 
$$
because $ \lambda_1 q^{-2n_2 } = \lambda_2 $. This finishes the proof of a).
Note that for this last implication we do not need $ l $ to be odd.
\medskip

We proceed to prove b).
We use the same principle for proving injectivity as in 
the proofs of theorem 1 and proposition 8 of [MR], although the combinatorial setup is different.

\medskip
Since $ {\cal G} $
generates $ b_n e $ as a right $ e b_n e $-module it induces a generating set 
of $ G \circ F M_n(\lambda) $ as a vector space
$$ {\cal M } := {\cal G } \otimes_{e b_n e } \seq_{n-2}^{n_1 -1 } \underline{12} $$

\medskip
We then have $ I := \varphi_{\lambda}({ \cal M }) = I_1 \cup I_2 $, where $ I_1 $ and $ I_2 $ are as above. 
Let us say that the elements of $ I_1 $ are of TL-type.
The elements of $ I $ are not independent: there are trivial relations between the TL-type elements as follows
$$
(Triv_1) \,\, \,  q^{-1} \,  w_1 12 w_2 \underline{12} w_3  -  w_1 21  w_2 \underline{12} w_3 = 
q^{-1} \,  w_1 \underline{12} w_2 12 w_3 -   w_1 \underline{12} w_2 21  w_3 
$$
for $ w_1 ,  w_2, w_3 $ words in $ 1 $ and $ 2 $, i.e. belonging to appropriate $ \seq_r^s $, 
% Let $ w_1 $ be of length $ i-1 $ and $ w_1 12 w_2 $ of length $ j-1$. Then the relations have preimages in $ \cal B $ of the 
% form 
% $$
% \begin{array}{c}
%  U_j \ldots U_{n-3}  U_{n-2} \otimes  ( q^{-1} \,  w_1 12 w_2  w_3 \underline{12} -  w_1 21  w_2  w_3 \underline{12} ) = \\ 
%  U_i U_{i+1} \ldots U_j \ldots U_{n-3}  U_{n-2} \otimes  ( q^{-1} \,  w_1 w_2 12 w_3 \underline{12} -  w_1   w_2 12  w_3 \underline{12} ) 
% \end{array}
% $$
% Under $ \pi $ these relations become the following trivial ones 
% $$
% \begin{array}{l}
%  U_i \, U_j U_{j+1} \ldots U_{n-3}  U_{n-2} \otimes_{e b_n e }  w_1 12 w_2 w_3 \underline{12}  = \\ 
%  U_i \,  U_j U_{j+1} \ldots U_{n-3}  U_{n-2} \otimes_{e b_n e }  w_1 12 w_2 w_3 \underline{12}  
% \end{array}
% $$
% Now, as in the [MR]-setup, using the relations $ \ast $ between the TL-type elements 
% as straightening rules, any element of $ I_2 $ can be put on the form $$ 22 \ldots 22 11 \ldots \underline{12} i_k i_{k+1} \ldots $$ i.e. with 
% no $ 12 $ appearing before the $ \underline{12}$. These elements are linearly independent.
% \medskip
\medskip

There are also certain trivial relations involving the first element $ U_{0, \ldots, n-1} :=  U_0 \, U_1 \cdots U_{n-1} $
of $ \cal G $ and the TL-elements. To handle these define first $ U_{0, \ldots, n-1}^{\lambda_1} := (U_0 + \lambda_1 ) \, U_1 \cdots U_{n-1} $
and replace then  $  U_{0, \ldots, n-1 }$ by 
$${\cal U }_{0, \ldots, n-1 }= (U_{n-1}+q) \, (U_{n-2}+q)  \cdots \,(U_{1} +q) U_{0, \ldots, n-1}^{\lambda_1 }\,  $$ 
in $ \cal G $. By this, $ \cal G $ remains a generating set of $ b_n $ as $ e b_n e $-module, since 
the expansion of ${\cal U }_{0, \ldots, n-1 }  $ gives 
$ U_{0, \ldots, n-1} $ 
plus a linear combination of the other elements of $ \cal G $ modulo $ e b_n e $.

\medskip

Now $U_0 + \lambda_1 = X $ and $ U_i = T_{i+1} -q $ and so we get 
$$ \varphi_{\lambda} ( {\cal U }_{0, \ldots, n-1 } \otimes_{e b_n e} i_1 i_2 \ldots i_{n-2} \underline{12}) =   S_{n-1} S_{n-2} \cdots S_2 \varpi \,\underline{ 12 }  i_1 i_2 \cdots i_{n-2}
 $$
Let us denote these elements by 
$\underline{1} i_1 i_2 \ldots  i_{n-2} \underline{2} $. They are 
$$\underline{1} i_1 i_2 \ldots  i_{n-2} \underline{2}
 := - \lambda_2 \,q^{n_2 -1 } 1 i_1 \ldots i_{n-2}2  + \lambda_1 \,
q^{n_1 -2} 2 i_1 \ldots i_{n-2} 1  $$

\medskip
The trivial relations 
between the $\, \underline{1} i_1 i_2 \ldots i_{n-2} \underline{2}\, $ and the TL-type elements
are then 
$$( Triv_2) \, \, \, q^{-1} \underline{1} w_1 12 w_2  \underline{2} \, - \, \underline{1} w_1 21  w_2 \underline{2} =
- \lambda_2 \,q^{n_2 -1 } 1 w_1 \underline{12} w_2  2 \, + \lambda_1 q^{n_1 -2} \, 2 w_1 \underline{12}  w_2 1 $$
where 
$ w_1 ,  w_2 $ are words in $ 1 $ and $ 2 $ belonging to appropriate $ \seq_r^s $, 

\medskip
To get a better understanding of these trivial relations we now consider $ w_1 \underline{12} w_2, \underline{1} w_3 \underline{2} $ as symbols and define 
$$ W_1 := \Span_k \{ w_1 \underline{12} w_2, \underline{1} w_3 \underline{2} \, | \, w_1 \in \seq_k^l, w_2 \in \seq_{n-k}^{l-n_1}, w_3 \in \seq_n^{n_1} \, \} $$ 
and $ W:= W_1/ \Span_k \{ R \, | \, R \in Triv_1 \cup Triv_2 \} $.
One checks on the relations that there is a linear map $ \psi_{\lambda}: \,  W \rightarrow G \circ F M_n(\lambda ) $ defined by
$$ \begin{array}{l} 
 w_1 \underline{12} w_2 \mapsto U_i U_{i+1} \ldots U_{n-1} \otimes_{e b_n e } w_1 w_2 \underline{12}, \\
\underline{1} w_3 \underline{2} \mapsto {\cal U }_{0, \ldots, n-1 } \otimes_{e b_n e } w_3 \underline{12} 
\end{array}
$$

Using the relations $ Triv_1 $ and $  Triv_2 $, it is straightforward to check that the elements 
$ 22 \ldots 11 \ldots 11\underline{12} i_k i_{k+1} \ldots i_n $ 
(with no $ 12  $ before the underline) and $ \underline{1} 222 \ldots 
111  \underline{2} $ generate $ W$.
We show that these elements 
map to a basis of $ M_n(\lambda) $
under $ \varphi_{\lambda} \circ  \psi_{\lambda} $ which implies that $ \varphi_{\lambda}  $ is injective.

\medskip
We have that
$$ \begin{array}{l}
\varphi_{\lambda} \circ  \psi_{\lambda} ( 22 \ldots 111 \underline{12} i_k  \ldots i_n ) = 22 \ldots 111 \underline{12} i_k  \ldots i_n \in M_n(\lambda)
\\ 
\varphi_{\lambda} \circ  \psi_{\lambda} (\underline{1} 222 \ldots 
111  \underline{2}) = \underline{1} 222 \ldots 
111  \underline{2} \in M_n(\lambda)
\end{array}
$$
The first kind of elements (of TL-type) were shown to be linearly independent in [MW1]. To show that 
$ \underline{1} 222 \ldots 
111  \underline{2} $ is independent of these, it is enough to show that it is nonzero modulo the TL-type elements. 
Calculating modulo the TL elements, we 
have $ 12 = q21 $ and so we find that $\underline{1} 222 \ldots 111  \underline{2}$ is equal to 
$$
\begin{array}{rl}
\underline{1} 2^{n_2 -1}  1^{n_1 -1}  \underline{2} & =  
- \lambda_2 \,q^{n_2 -1 } \, 1 2^{n_2 -1}  1^{n_1 -1} 2  + \lambda_1 \,
q^{n_1 -2} \, 2^{n_2 }  1^{n_1 }  \\ & 
= (-\lambda_2  q^{2n_2 + n_1-2}    + \lambda_1 q^{n_1-2}   ) 2^{n_2+1} 1^{n_1 +1}   
\end{array}
$$
By the assumption of the lemma this is nonzero since $ \lambda_1 / \lambda_2 = q^{ 2 m } $.

\medskip
Finally the other implication of b) follows also from the last calculation since $ \psi_{\lambda} $ is surjective. We have proved the lemma.
\end{pf*}

\medskip
A consequence of the lemma is that $ M_n(\lambda) $ is not isomorphic to $ \Delta_n ( \lambda) $ in general. 
Moreover, we shall later in section 5 explain how the above proof can be used to deduce that 
$ M_n(\lambda) $ is also not isomorphic to $ \nabla_n ( \lambda) $ in general.

\medskip
On the other hand, we now prove by induction that $ M_n( \lambda) $ and $ \Delta_n ( \lambda) $ 
are equal in the Grothendieck group of 
$ b_n $-modules. The next lemma is the induction basis. 

\begin{lemma}{\label{inductionbasis}} 
For $ n \geq 1 $ we have the following isomorphisms of $ b_n $-modules
$$ \begin{array}{lll}
a) \,\,\, M_n(n) \cong \Delta_n(n)  & 
b) \,\,\,  M_n(-n) \cong \Delta_n(-n) & 
c) \,\,\, M_2(0) \cong \Delta_2(0)  
\end{array} $$
\end{lemma} 

\begin{pf*}{\it Proof} 
The parts a) and b) of the lemma are easy to check since all the involved $ b_n $-modules 
are one dimensional and have trivial $ U_i $ actions for $ i \geq 1$. One then just needs verify that $ U_0 = X - \lambda_1 $ 
acts the right way.

\medskip

In order to prove part c) we first get  
for $ n=2$ by direct calculations that the matrices of $ U_1 $ and 
$ X $ with respect to the basis $ \{ 12, 21 \} $ of $M_2(0)$ are
given by
$$ 
\begin{array}{rl}
 U_1  =   \begin{pmatrix}
   -q^{-1} & 1 \\
   1 & -q \\
\end{pmatrix},  & 
  X  =   \begin{pmatrix}
   \lambda_1 & 0 \\
   -\lambda_1 (q-q^{-1}) & \lambda_2 \\
\end{pmatrix},
\end{array} 
$$
and hence the matrix of $ U_0 = X - \lambda_1 $ is 
$$  U_0  =   \begin{pmatrix}
   0 & 0 \\
   -\lambda_1 (q-q^{-1}) & -[m] \\
\end{pmatrix}$$
since $ [m] = \lambda_1 - \lambda_2$. 
The ket basis of $ \Delta_2(0) $, see [MW1], modulo multiplication by nonzero scalars, is given by $ \{ \cup, U_0 \, \cup  \}$.
Define $ \varphi $ by 
$$ \varphi: \underline{12} = q^{-1}12 -21 \mapsto \cup, \,\,\,\,\, \, \, \, U_0 \,\underline{12} \mapsto U_0 \,\cup $$
This is the desired $ b_n $-isomorphism provided that $ U_0 \, \underline{12}  $ is nonzero 
and is an eigenvector of $ U_0 $ with eigenvalue $ -[m] $.
But by the above $$ U_0 \, \underline{12} =  q^{-1} (-\lambda_1 (q-q^{-1}) + q [m] ) 21 $$
The coefficient is nonzero iff $ \lambda_1 (q-q^{-1}) \not= q [m] $, which by $ \lambda_1 = \frac{q^m}{q-q^{-1} }$ is 
equivalent to $ q^{2m} \not= q^2 $, which holds by the assumptions on $q $ 
given in the beginning of section 3. 
But then 
$ \underline{12} $ is automatically an eigenvector for $ U_0 $ of the right eigenvalue.
\end{pf*}

\begin{Thm}{\label{Grothendieck}}
Assume that $ n \geq 1 $. Then $[ \Delta_n( \lambda): L_n(\mu) ] 
=  [ M_n( \lambda): L_n(\mu) ] $ for all $ \lambda,\mu \in \Lambda_n  $.
\end{Thm}
\begin{pf*}{\it Proof} We prove the theorem by induction on $ n $. The induction basis $ n= 1 $ and $ n= 2 $ is provided by the above lemma. 
We assume the theorem to hold for 
all $ n^{\prime} $ strictly smaller than $ n $ and prove it for $ n  $.  
Recall once again that the simple $b_n$-modules $ L_n( \mu) $ satisfy that 
$$
\begin{array}{l}
F L_n(\mu) \cong \left\{ \begin{array}{ll} 
L_{n-2}(\mu) & \mbox{if} \,\,\,   \mu \in \Lambda_{n}  \setminus \{\pm n\}\\
0  & \mbox{otherwise}  
\end{array}
\right. 
\end{array}
$$
By induction, exactness of $ F $, the category theoretical 
property for $ \Delta_n(\lambda) $ stated in 
(\ref{categorial}) and 
Theorem (\ref{firststep}), we then 
get for $ \mu \in \Lambda_{n}  \setminus \{\pm n\} $ that
$$ 
\begin{array}{l}
\, [ \Delta_n(\lambda): L_n(\mu) ] = [ F \Delta_n(\lambda) :  F L_n(\mu) ] = [  \Delta_{n-2} 
(\lambda): L_{n-2}(\mu) ] = 
\\ \,  [  M_{n-2} (\lambda) : L_{n-2}(\mu) ] =  [ F  M_{n} (\lambda) :F  L_{n}(\mu) ] =  [   M_{n} (\lambda) :  L_{n}(\mu) ] $$
\end{array} 
$$
and we need now only prove $ [ \Delta_n(\lambda) : L_n(\pm n ) ] = [ M_n(\lambda): L_n(\pm n ) ] $.

\medskip 
But $ X $ acts semisimply in any $ b_n $-module and so we obtain the following $ {\C}[ X] $-module
decompositions
$$ \Delta_n(\lambda ) = \bigoplus_{\mu \in \Lambda_n }  L_n(\mu)^{ d_{\lambda \mu} }, \,\,\,\,\, 
M_n(\lambda ) = \bigoplus_{\mu \in \Lambda_n }  L_n(\mu)^{ e_{\lambda \mu} } $$
where $ d_{\lambda \mu} = [ \Delta_n(\lambda): L_n(\mu) ] $ and $ e_{\lambda \mu} = [ M_n(\lambda): L_n(\mu) ] $.
On the other hand, the only possible eigenvalues for $ X $ are $ \lambda_1 $ and $ \lambda_2 $ 
and we just saw that 
$ d_{\lambda \mu } = e_{\lambda \mu } $ for $ \mu \in \Lambda_n \setminus \{ \pm n \} $.
Hence it is enough to show that 
$ \Delta_n(\lambda) $ and $  M_n(\lambda ) $ are isomorphic  $ {\C}[ X] $-modules 
to deduce 
$ d_{\lambda \mu } = e_{\lambda \mu } $ for the remaining $ \mu \in \Lambda_n $
and so finish the proof.
Indeed $ L_n(n) $ and $ L_n(- n) $ are 
both one dimensional, generated by eigenvectors for $ X $ of eigenvalues $ \lambda_1 $ and 
$ \lambda_2 $ respectively (recall $ \lambda_1 \not= \lambda_2 $ by our assumptions).

\medskip
Now $ \Delta_n(\lambda) \cong  M_n(\lambda ) $ as 
$ {\C}[ X] $-modules
if and only if the eigenspace multiplicities with respect to $ X $ are equal, so 
we show that this is the case.

\medskip
For this we observe that the Bratteli diagram or Pascal triangle 
of restriction rules from $ b_n $ to $ b_{n-1} $ given in [MW1] can be used to determine the eigenvalues of $ X $ on $ \Delta_n(\lambda) $ 
in the following way: A diagram of the diagram basis of 
$ \Delta_n(\lambda) $ 
is an eigenvector for $ X = U_0 + \lambda_1 $ of eigenvalue $ \lambda_2 $ iff 
its first line is marked with a filled blob.
This induces the 
following Pascal triangle pattern of multiplicities of the eigenvalue $ \lambda_2 $.
$$
\begin{array}{llllllllll}
n=1 \, \, \, \, \, \, & & & &1 & &0 & & & \\
n=2 \, \, \, \, \, \, & & &1&  &1 & & 0& & \\
n=3 \, \, \, \, \, \, & &1& & 2&  &1& &0 &\\
n=4 \, \, \, \, \, \, &1 & &3 && 3 &&1 & &0
\end{array}
$$
For example, the first number $3$ says that 
$ \Delta_4(-2) $ has $3$ diagrams with first line marked and hence $ \lambda_2 $ has multiplicity 3 
in $ \Delta_4(-2) $.

\medskip
We must compare this pattern with the $ \lambda_2 $-multiplicity of $ X $ in $ M_n(\lambda) $.
We have with the usual notation $ \lambda = n_1 - n_2 $ 
a basis of $ M_n(\lambda) $ consisting of
$ B := \seq_n^{n_1} $. Define $ B_1 $ as the sequences from $ \seq_n^{n_1} $ 
that begin with a $ 1 $ and 
$ B_2 $ as $ \seq_n^{n_1} \,  \setminus \,  B_1 $.  Put an order on $ B $ such that 
the elements of $ B_2 $ come before the elements of $ B_1 $. Then by lemma \ref{basiclemma} 
the action of $ X $ is upper triagonal with $  \lambda_2 $ in the first $ | B_2 | $ diagonal 
elements and with $ \lambda_1 $ in the last $ | B_1 | $ diagonal elements. Hence the 
$ \lambda_2 $-multiplicity of $ X $ is $  | B_2 | $. But the numbers $ B_2 $ satisfy the 
same Pascal triangle recursion as the above. The theorem is proved.

\end{pf*}
\section{\bf Specht modules and duality}
In this section we shall relate the results of the previous sections to 
the $ H_k(n,2)$-module $ \tilde{S}^{\lambda} $
introduced in [DJM] for bipartitions $ \lambda =(\tau, \mu) $ 
of $ n $.
The module $ \tilde{S}^{\lambda} $ is a cell module for a certain cellular structure 
on $ H_k(n,2)$. Following 
modern terminology as used in for example [Ma], we shall therefore denote it the 
{\it Specht module} for $ H_k(n,2)$, 
although it is rather an analogue of the dual Specht module, 
and for $ \lambda =(\tau, \mu) $ 
we shall 
accordingly use the notation $ S(\lambda) $ or $ S( \tau, \mu) $
for it.	 If $ \lambda = ( (n_1), (n_2))$
is a two-line bipartition of $ n $, that is $ n_1, n_2 \geq 0 $  
such that $ n_1 + n_2 = n $, we 
shall also write $  S(n_1, n_2) $ for $ S(\lambda) $. Similarly, if 
$ \lambda = ( (1^{n_1}), (1^{n_2}))$ is a two-column bipartition, we shall write 
$  S(1^{n_1}, 1^{n_2}) $ for $ S(\lambda) $.

\medskip
In this section we show that the Specht module $  S(n_1, n_2) $ as 
well as its contragredient dual $ S(n_1, n_2)^{\circledast}$
are modules for $ b_n$. 
We moreover establish a $ b_n$-isomorphism between $  S(n_1, n_2)^{\circledast} $
and $M_n(\lambda)$ 
where $ \lambda = n_1 -n_2 $. Finally, we prove an analogue of lemma \ref{relate} 
for $ M_n( \lambda)^{\circledast} $ and as a consequence we get that, somewhat surprisingly, neither 
of the modules $  S(n_1, n_2), S(n_1, n_2)^{\circledast}, 
M_n( \lambda), M_n( \lambda)^{\circledast} $ is the pullback of the standard module
$ \Delta_n(\lambda) $ for $ b_n$ in general.

\medskip
On the other hand, the pullback of the simple $ b_n $-module
$ L_n(\lambda) $ to 
$ H_k(n,2 ) $ certainly {\it is} a simple $ H_k(n,2 ) $-module. 
Thus, the statements of the previous paragraph are 
apparently not compatible with the statement of 
theorem \ref{Grothendieck} on equality in the Grothendieck groups, 
since the dominance order on bipartitions does not induce the quasi-hereditary order
$ \prec $ on $ \Lambda_n $. But note that 
the bipartitions $ (\tau , \mu ) = ( (n_1), (n_2) ) $ are only Kleshchev (= restricted) 
in 'small' cases and therefore, apart from these 
small cases, $ L_n(\lambda) $ 
is not the simple module associated with the bipartition $( (n_1), (n_2) ) $ 
when viewed as $ H_k(n,2) $-module, see [AJ].
In fact, it would be interesting 
to know which is the Kleshchev bipartition corresponding to $ L_n(\lambda) $.
(In the recent preprint [RH] we have solved this problem).

\medskip
Let us now recall the combinatorial 
description of the permutation module $ M_{H}(\tau , \mu) $ 
and the Specht module $ S(\tau , \mu) $ for $ H_k(n,2)$ given in 
[DJM] and 
[DJMa]. Since these references use right modules rather than left modules and since they 
moreover use a slightly different presentation of $ H_k(n,2)$, 
the following formulas vary slightly from theirs.

\medskip
Let $ (\tau , \mu) $ be a bipartition of $n$. Then 
a $ (\tau , \mu) $-bitableau $ t $ is a pair $  (t^1,t^2) $ where $ t^1 $ is a $ \tau $-tableau and  $ t^2 $ is a $ \mu$-tableau 
and where tableaux means fillings with the numbers $ I_n = \{ \pm 1, \pm 2, \ldots ,\pm n \} $, where either $  i $ or $ -i $ occurs exactly once. 
Two $ (\tau , \mu) $-bitableaux 
$  (t^1,t^2) $ and $  (s^1,s^2) $ are
said to be row equivalent if the tableaux obtained by taking absolute values in $ t^1 $ and $ s^1 $ are row equivalent in the usual sense, and if $ t^2 $ and $ s^2 $ are row equivalent.
The equivalent class of the bitableau $ t $ is called a tabloid 
and is written $ \{t\} $.

\medskip
The permutation module $ M_H(\tau , \mu ) $ 
for $ H_k(n, 2) $ is now  
$$ M_H(\tau , \mu ) := \Span_k \{ \{ t_1,t_2 \} \, | \, (t_1,t_2) \mbox{ is a 
row standard } (\tau , \mu )
\mbox{-bitableaux } \} $$
where the action can be read off from the lemmas 3.9, 3.10 and 3.11 of [DJMa].

\medskip
The Specht module $ S_H(\tau , \mu) $ is now the quotient $ M_H(\tau , \mu)/ N_H(\tau , \mu) $ 
for $ N_H(\tau , \mu) $ a certain submodule 
of $ M_H(\tau , \mu) $. The standard tabloids induce a basis for $ S(\tau , \mu) $ 
$$ [ t_1,t_2] :=  \{ t_1,t_2 \} +  N_H(\tau , \mu) $$ 
where standard means that all entries are positive, 
and that each component is row standard and column standard. 

\medskip 
We shall be especially concerned with the case of two-line bipartitions
$ (\tau , \mu ) = ((n_1), (n_2 ) )$. In that case, standard bitableaux are just row standard 
tableaux with positive entries and so the 
formulas for the action of $ H_k(n,2) $ on $ M_H (\tau, \mu) $ 
induce the following formulas for the action on 
$ [t]= [t_1, t_2] \in S(\tau , \mu) $ 
\begin{equation}{\label{action_rules}} g_i [t]= 
  \left\{ \begin{array}{ll}
  \sigma_i [t ]  & \mbox{if} \,\,\,   (i \in t^1, i+1 \in t^2)       \\
  \sigma_i [t ] + (q-q^{-1}) [  t ] & \mbox{if} \,\,\,   (i+1 \in t^1, i \in t^2)        \\
 q [  t ]  & \mbox{if} \,\,\,   (i,i+1 \in t^1 ) \, \, \mbox{or} \,\,\,    
(i, i+1 \in t^2)       \\
\end{array}
\right.
\end{equation}
where the transposition $ \sigma_i = (i,i+1) $ acts by permuting the entries. 
The action of $ X $ can only partially be made explicit. 
We consider first the action of $ X_i $.  Let $ t^{\tau , \mu } $ be the $ (\tau, \mu)$-bitableau with 
$ \{ 1, \ldots , n \} $ positioned increasingly from left to right.
For example, in the case $ n_1 = 5 , n_2 = 6 $ we have 
$$
t^{\tau  , \mu }=(\,  {\small{\tableau[scY]{1,2, 3, 4, 5 |}}}\, , {\small{\tableau[scY]{6,7,8,9,10,11 |}}} \, ) 
$$

Then by [DJM] we have 
$$ X_i [ t^{\tau , \mu } ] = 
\left\{ 
\begin{array}{ll}
\lambda_1 q^{2(i-1)} [ t^{\tau , \mu } ]  & \mbox{if} \, \, \,  
i = 1, \ldots , n_1 \\
\lambda_2 q^{2 (i-n_1-1)} [ t^{\tau , \mu } ]   & \mbox{if} \, \, \,  
i = n_1+1, \ldots , n 
\end{array} \right. 
$$
To get the action on the other standard tableaux, one has to use 
the commutation rules of $ H_n(n,2) $. 
This implicit description is enough to prove the following theorem.  Although it is a main  
philosophical idea of [MW], a formal proof was not given.
\begin{Thm}
$ S(\tau, \mu) $ is a module for $ b_n $ when $ (\tau , \mu ) = ((n_1), (n_2 ) )$.
%The pull-back of the standard module $ \Delta_n( 2k -n ) $ is the Specht module $ S((k), (n-k)) $.
\end{Thm}
\begin{pf*}{\it Proof}
By the isomorphism theorem (\ref{MW-iso}) we must verify that 
\begin{equation}{\label{check}} (X_1 X_2 - \lambda_1 \lambda_2 ) (g_1 -q ) =0
\end{equation}
in $ \End_k (S(n_1, n_2)) $. 
Let therefore $ [ t ] = [ t_1, t_2 ]  $ be the class of a standard bitableau for the 
bipartition $ ((n_1), (n_2)) $. If $ 1, 2 $ both
belong to $ t_1 $ or $ t_2 $ the statement is clear by (\ref{action_rules}).
Using (\ref{action_rules}) once again,
we have that 
$$ S( n_1, n_2) =  \ker(g_1 -q) + \Span_k\{ \, [ t_1, t_2 ] \, | \, 1 \in t_1 , \, 2 \in t_2   \}  $$
and we are left with the case $ 1 \in t_1 $, $ 2 \in t_2 $.
But then we can find $ w = \sigma_{i_1} \ldots \sigma_{i_r} \in \langle \sigma_i \, | \, i = 2, \ldots,  n-1 \rangle $ 
such that $ w \, t^{\tau , \mu }  = (t_1, t_2 ) $ and so we have 
$ X_1 [ t_1, t_2 ] = \lambda_1 [ t_1, t_2 ] $ since $ X= X_1 $ commutes 
with all $ g_2, \ldots , g_{n-1} $. 
 
\medskip
We then consider the action of $ X_2 $ on $ [ t_1, t_2 ] $. Let $  t^{12} $ be the bitableau with $ 1 \in t^1 $, $ 2 \in t^2 $ and the other 
entries increasing from left to right. For example, if $ n_1 = 5 $ and $ n_2 = 6 $, it is 
$$
t^{12}=(\,  {\small{\tableau[scY]{1,3,4,5,6 |}}}\, , {\small{\tableau[scY]{2,7,8,9,10,11 |}}} \, ) 
$$
Then any $ t = (t_1, t_2 ) $ with $ 1 \in t_1 $ and $ 2 \in t_2 $ is of the form $ t = w \, t^{12} $ 
where 
$ w = \sigma_{i_1} \ldots \sigma_{i_r}  \in \langle \sigma_i \, | \, i = 3, \ldots,  n-1 \rangle $. 
We claim that $ X_2 [ t^{12} ] = \lambda_2 [ t^{12} ]$ modulo a linear combination of 
elements $ [ (t^1, t^2 ) ] $ all satisfying $ 1,2 \in t^1 $. Believing this, 
we would also get that $ X_2 [ t ] = \lambda_2 [ t ]$ modulo a similar linear combination 
of elements $ [ (t^1, t^2 ) ] $, since $ X_2 = g_1 X g_1 $ and $ g_i $ 
commute for $ i = 3, \ldots, n $. From this (\ref{check}) would follow.

\medskip 
To prove the claim for $ t^{12} $ we first use (\ref{action_rules}) to write
$$ g_{2} \,g_3 \ldots    \, g_{n_1 -1 } \,  g_{n_1} \,  \{ t^{\tau  \mu } \} = \{ t^{12} \} $$ 
Since $ X_{n_1 +1}^{-1} \{ t^{\tau  \mu } \} = \lambda_2^{-1} \{ t^{\tau  \mu } \} $ and 
$ X_{n_1+1} =    \,g_{n_1} \ldots    \, \,  g_1     X_1   g_{1} \ldots    \,  g_{n_1} $ we deduce that 
$$ X_2 \{ t^{12}\} = \lambda_2 \,  g_2^{-1} \, \ldots g_{n_1}^{-1} \,   \{ t^{ \tau  \mu } \} $$  
The claim now follows.
\end{pf*}

\medskip
Recall that the contragredient dual $ M^{\circledast}$ of an $H_k (n, 2)   $-module $M $
is the linear dual $ \Hom_k (M, k ) $ 
equipped with the $ H_k (n ,  2) $ action $ (h f)(m):= f(h^* m) $ 
for 
$ * $ the antiinvolution of $ H_k (n,  2) $ given by $ g_i^* := g_i $ and $ X^*:= X $.

\medskip
Let $ H_k^{\prime} (n ,  2) $ be the Ariki-Koike algebra $ H_k ( -q^{-1}, \lambda_2, \lambda_1 ) $.
There is a $ k$-algebra isomorphism $ \theta:  H_k (n ,  2) \rightarrow H_k^{\prime} (n ,  2)  $ 
given by 
$$ X \mapsto X, \,\,\, \, \, \,  \, \, \, \, \, g_i \mapsto g_i $$
Following [Ma] and [F], we define $ S^{\prime} (\tau , \mu ) $ 
as the pullback under  $ \theta $ of the Specht module  $ S (\tau , \mu ) $ for
$ H_k^{\prime} (n ,  2) $. Now Mathas proved in [Ma] the following result.

\begin{Thm}
As $ H_k(n,2) $-modules we have $S( \tau , \mu )^{\circledast}  \cong S^{\prime} ( \mu^{\prime}, \tau ^{\prime} ) $
where $ \tau ^{\prime} $ and $ \mu^{\prime} $ are the usual conjugate partitions of $ \tau  $ and $ \mu$.
\end{Thm}

In the case $ (\tau , \mu ) = ((n_1), (n_2) ) $, 
the isomorphism of the theorem will also be an isomorphism of $ b_n $-modules, since $ * $ induces the usual antiinvolution $ * $ of $ b_n$ that appears 
in the definition of contragredient duality in $ b_n$-mod. Specially, 
$ S^{\prime} ( 1^{n_2} , 1^{n_1} ) $ will be a $ b_n $-module as well.

\medskip
The standard basis for $ S(\mu^{\prime} , \tau ^{\prime}) =   S^{\prime} ( 1^{n_2} , 1^{n_1} ) $ 
consists of the classes of bitableaux $ t = (t_1, t_2 ) $ of the 
bipartition $ ((1^{n_2}) , (1^{n_1})) $. 
We get for $ g_i $ the same 
action rules as before:
\begin{equation} g_i [t]= 
  \left\{ \begin{array}{ll}
  \sigma_i [t ]  & \mbox{if} \,\,\,   (i \in t^1, i+1 \in t^2)       \\
  \sigma_i [t ] + (q-q^{-1}) [  t ] & \mbox{if} \,\,\,   (i+1 \in t^1, i \in t^2)        \\
 q [  t ]  & \mbox{if} \,\,\,   (i,i+1 \in t^1 ) \, \, \mbox{or} \,\,\,    
(i, i+1 \in t^2)       \\
\end{array}
\right.
\end{equation}
As before, we have a special standard bitableau $ t^{\mu^{\prime} , \tau ^{\prime}} $, 
this time with the numbers $ 1, \ldots , n $ filled in increasingly 
first down the first column, then down the second column. 
The action of $ X_i $ on this $[t^{  \mu^{\prime} , \tau ^{\prime} }] $ is given by 
$$ X_i [ t^{  \mu^{\prime} , \tau ^{\prime}         } ] = 
\left\{ 
\begin{array}{ll}
\lambda_2 q^{2(i-1)} [ t^{  \mu^{\prime} , \tau ^{\prime} }  ] & \mbox{if} \, \, \,  
i = 1, \ldots ,n_2 \\
\lambda_1 q^{2 (i-n_2 -1 )} [ t^{  \mu^{\prime} , \tau ^{\prime} } ]   & \mbox{if} \, \, \,  
i = n_2+1 , \ldots , n 
\end{array} \right. 
$$

We are now in position to prove the following result
\begin{Thm}{\label{duality}}
Let as before $ \lambda = n_1 -n_2 $.
Then there is an isomorphism of $ b_n $-modules 
$ M_n (\lambda) \cong S( n_1, n_2 )^{\circledast}  $.
\end{Thm}
\begin{pf*}{\it Proof}
We had by Mathas's theorem that 
$ S( n_1, n_2 )^{\circledast}  \cong S^{\prime} ( 1^{n_2 }, 1^{n_1} ) $. 
We then define a linear 
map $ \varphi : S^{\prime} ( 1^{n_2 }, 1^{n_1} ) \rightarrow M_n (\lambda) $ 
by $$ \varphi ( [ t_1 , t_2 ] )= i_1 i_2 \ldots i_n   \, \, \mbox{ where  }   
i_j = 1 \, \,  \mbox{ iff } j \in t_2$$
It is easily checked that $ \varphi $ is linear with respect to $ g_i $. On the other hand, we have that 
 $ \varphi(  t^{  \mu^{\prime} , \tau^{\prime} }) =  
2^{n_2} 1^{n_1} $. 
Using the next lemma we see that $ X_i $ acts through the same constant on 
$   [ t^{  \mu^{\prime} , \tau^{\prime} } ]  $ and as on 
$ 2^{n_2} 1^{n_1} $. This is enough to complete the proof by the commutation rules for $ H_k(n,2) $.
\end{pf*}

\begin{lemma}{\label{X-action}} Let $ w = 2^{n_2} 1^{n_1} \in M_n(\lambda) $. Then 
$$ X_i w  = 
\left\{ 
\begin{array}{ll}
\lambda_2 q^{2(i-1)} \, w   & \mbox{if} \, \, \,  
i = 1, \ldots ,n_2  \\
\lambda_1 q^{2 (i-n_2 -1 )} \, w    & \mbox{if} \, \, \,  
i = n_2+1 , \ldots , n 
\end{array} \right. 
$$
\end{lemma}
\begin{pf*}{\it Proof}
By lemma \ref{basiclemma} the action of $ X =  X_1 $ on $ w $ is multiplication by $ \lambda_2  $, hence the action of $ X_2 = T_2 X_1 T_2 $ is 
multiplication by $ q^2 \lambda_2 $ and so on until we reach $ X_{n_2} $. 

\medskip
To calculate the action of $ X_{n_2+ 1} $ we write
$$ 
\begin{array}{r} X_{n_2 +1} = 
T_{n_2 + 2}^{-1}  \ldots  T_n^{-1} S_n \ldots S_2 \varpi T_2 \ldots T_{n_2 + 1} 
\end{array}
$$
and so 
$$
\begin{array}{l}
X_{n_2+1} w  = 
T_{n_2 + 2}^{-1}  \ldots  T_n^{-1} S_n \ldots S_2 \varpi T_2 \ldots T_{n_2 + 1} 2^{n_2} 1^{n_1} = \\
\lambda_1 T_{n_2 + 2}^{-1}  \ldots  T_n^{-1} S_n \ldots S_2  1 2^{n_2} 1^{n_1-1} =
q^{n_1 -1} \lambda_1 T_{n_2 + 2}^{-1}  \ldots  T_n^{-1} 2^{n_2 } 1^{n_1} = \\ 
\lambda_1 2^{n_2} 1^{n_1} = \lambda_1 w 
\end{array}
$$
and the action is multiplication by $ \lambda_1$.
This implies that $ X_{n_2 +2 } $ acts by $ \lambda_1 q^2 $ and so on. 
\end{pf*}

We can now finally prove the result alluded to in the previous section.
\begin{Cor} 
Let $ n \geq 3 $ and suppose that $ q $ is an $l$th primitive root of unity, where
$ l $ is odd.
Suppose $ \lambda \in \Lambda_{n} \setminus \{\pm n \} $. Then the 
adjointness 
map $ \psi_{\lambda}: G \circ F M_n(\lambda)^{\circledast} \rightarrow M_n(\lambda)^{\circledast} $ is 
an isomorphism iff $ n_1 = m \mbox{ mod } l $. 
\end{Cor}
\begin{pf*}{\it Proof}
By the actions rules given above and theorem \ref{duality} the actions on $ M_n(\lambda)^{\circledast} $ and $ M_n(\lambda) $ 
are the same, except that $ \lambda_1 $ and $ \lambda_2 $ are interchanged as are $ n_1 $ and $ n_2$. We then repeat the argument of Lemma \ref{relate} and get that $ \varphi_{\lambda} $ is 
an isomorphism iff $ \lambda_2 / \lambda_1 = (-q)^{-2 n_1} $, which is equivalent to 
$ n_1 = m \mbox{ mod } l $ as claimed.
\end{pf*}

Combining the corollary with lemma \ref{relate}
we deduce that neither $ M_n (\lambda) $ nor $ M_n(\lambda)^{\circledast} $ 
is the standard module $ \Delta_n(\lambda) $ for $ b_n$ in general. 
And then, combining this with 
the above theorem, we get the same statement for the Specht module $ S_n(n_1, n_2 ) $ and for 
$ S_n(n_1, n_2 )^{\circledast} $.

\section{\bf Alcove geometry}
We already saw that although $ M_n(\lambda) $ does not identify with the standard 
module $ \Delta_n(\lambda) $ for $ b_n $ in general, the two modules 	still 
have many features in common. In this section we shall further pursue this point, by considering the 
behavior of the restriction functor $ \res^{b_n}_{b_{n-1}} $ from $ b_n $-mod to $ b_{n-1} $-mod on $ M_n(\lambda)$.

\medskip

It is known from [MW1] that the representation theory of $ b_n $ is governed by an alcove geometry on $ \mathbb Z $ where 
$ l $ determines the alcove length and $ m $ the position of the fundamental alcove. There is a linkage principle 
and the decomposition numbers are given by Kazhdan-Lusztig polynomials for the affine Weyl group
$ \cal W $ corresponding to $ {\mathfrak sl}_2  $.

\medskip

Let us now set up some exact sequences that arise from 
restriction from $ b_n $-mod to $ b_{n-1} $-mod.
Let $ \lambda \in \Lambda_{n} \setminus \{ \pm n \}  $. 
As a $\TL_{n-1}$-module the restricted module 
$ \res^{b_n}_{b_{n-1}} M_n(\lambda ) $ is isomorphic to the direct sum $$ M_{n-1}(\lambda +1 ) \oplus
M_{n-1}(\lambda -1 ) $$
This is however not automatically the case when 
$ \res^{b_n}_{b_{n-1}} M_n(\lambda ) $ is considered as a $b_{n-1} $-module since $ X $ acts differently as 
element of $ b_n $ and of $ b_{n-1} $. But the following statement always holds.

\begin{lemma}{\label{sequences}} 
Assume $ \lambda \in \Lambda_{n} \setminus \{\pm n \} $. 
Then there is a short exact sequence 
of $ b_{n-1}$-modules 
$$ 0 \rightarrow M_{n-1}(\lambda - 1 ) \rightarrow 
 \res^{b_n}_{b_{n-1}} M_n(\lambda )  \rightarrow  M_{n-1}(\lambda + 1 ) \rightarrow 0 $$
\end{lemma}
\begin{pf*}{\it Proof}
We identify $ M_{n-1}(\lambda - 1 ) $ with the span of the sequences of the 
form $ v_{1} v_{2} \cdots v_{n-1} 1 $. Since
for all $ x \in \seq_{n-2} $ we have that $ T_n^{-1} S_n ( x 11) = x 11 $ and
$$ T_n^{-1} S_n ( x 21 ) = T_n^{-1} ( x 12  ) = x21 $$ 
we get that $ M_{n-1}(\lambda - 1 ) $ in this way is 
a $ b_{n-1} $-submodule of $ \res^{b_n}_{b_{n-1}} M_n(\lambda )$.

\medskip
The quotient 
of $ \res^{b_n}_{b_{n-1}} M_n(\lambda )$ 
by $ M_{n-1}(\lambda - 1 ) $ is now generated by the classes of the 
sequences that end in 2.
It can be identified with $ M_{n-1}(\lambda + 1 ) $ since 
for $ x \in \seq_{n-2} $ we have  $ T_n^{-1} S_n ( x 22) = x 22 $ and
$$ T_n^{-1} S_n ( x 12 ) = T_n^{-1} ( x 21 ) = x 12  \mod  M_{n-1}(\lambda - 1 ) $$
The lemma now follows.
\end{pf*}
One observes that these
sequences are very similar to the 
sequences for $ \res^{b_n}_{b_{n-1}} \Delta_n(\lambda )$ 
given in lemma 4.5 of [MW1]. The only difference is that in [MW1] 
the appearances of $ \lambda -1 $ and $ \lambda +1 $ are interchanged 
when $ \lambda $ is negative. 
But $ M_n(\lambda) $ is not the pullback of $ \Delta_n(\lambda) $, as we already pointed out
several times, 
and it seems to be a difficult task to compare the two 
systems of exact sequences. 

\medskip

We finish the paper by showing that the sequences of the 
lemma are split when $ \lambda $ is not a wall of the alcove geometry. 
This result could also have been obtained using theorem \ref{Grothendieck} and 
the linkage principle
for $ b_n $-mod, but we here 
deduce it from the machinery we have set up.
We use central elements.

\medskip
It is known,
see for example the appendix of [MW], that the symmetric polynomials in the $ X_i $ are central elements of $ H(n,2)$ and hence also of $ b_n $.
We consider $ z := X_1 X_2 \ldots X_n $ as an element of the centre $ Z(b_n ) $ of $ b_n$ and work out the 
action of it on $ M_n(\lambda) $.

\begin{lemma} Recall that $ \lambda = n_1-n_2$. Then the action of $ z$ on $ M_n(\lambda) $ is diagonal, given by the constant 
$$ \lambda_1^{n_1} \lambda_2^{n_2} q^{n_1 (n_1 -1) } q^{n_2 (n_2 -1) } $$
\end{lemma}

\begin{pf*}{\it Proof}
As a $ b_n $-module $ M(\lambda ) $ is generated by $ 2^{n_2} 1^{n_1} $. Since $ z$ is central, it is therefore enough to 
prove the assertion on that element.
Recall that the $ X_i $ commute.
By lemma \ref{X-action}
we find that $ X_1 X_2 \ldots X_{n_2} $ acts by 
$$ \lambda_2^{n_2} q^{0+2 + 4 + \ldots 2(n_2-1) } =
\lambda_2^{n_2} q^{n_2(n_2-1) }$$
Once again by lemma \ref{X-action}, we have 
that $ X_{n_2 +1 } \ldots X_{n} $ acts by
$$ \lambda_1^{n_1} q^{0+2 + 4 + \ldots 2(n_1-1) } =
\lambda_1^{n_1} q^{n_1(n_1-1) }$$
The lemma now follows by combining.
\end{pf*}

We can now prove the promised splitting. 
\begin{Thm}
Assuming $ \lambda \not= -m \! \mod  l $,
the exact sequences from lemma \ref{sequences} are split.
\end{Thm}
\begin{pf*}{\it Proof}
If the sequence were nonsplit, any preimage 
in 
$ \res^{b_n}_{b_{n-1}} M_n(\lambda ) $ of  
the $ M_n(\lambda+1) $ generator $ w =2^{n_2} 1^{n_1}  $ 
would generate a submodule $ M \subset \res^{b_n}_{b_{n-1}} M_n(\lambda ) $ nonisomorphic to $ M_n(\lambda+1) $.
Moreover $ M $ would map surjectively onto $ M_n(\lambda+1) $ and would have a composition factor in common with $ M_{n-1}( \lambda -1 ) $.
But then $ z $ would act through the same constant on $ M_n(\lambda+1) $ and $ M_n(\lambda-1) $.

\medskip
Let $ \lambda = n_1 - n_2$. The action of $ z $ on $ M_{n-1}( \lambda-1 ) $ is  
$$ \lambda_1^{n_1-1} \lambda_2^{n_2} q^{(n_1-1) (n_1 -2) } q^{n_2 (n_2 -1) } $$
and the action of $ z $ on $ M_{n-1}( \lambda+1 ) $ is 
$$ \lambda_1^{n_1} \lambda_2^{n_2-1} q^{n_1 (n_1 -1) } q^{(n_2-1) (n_2 -2) } $$
Equating, we get $$ \lambda_2 q^{2(n_2-1)} = \lambda_1 q^{2(n_1-1)} $$ 
which implies that $ \frac{ \lambda_1}{\lambda_2} = q^{2m} = q^{ 2(n_2 -n_1)} $ and the theorem follows.

\end{pf*}

\end{document}